\documentclass[MSNbibl,number,citesort,seceqn,dvips]{arxbj}
\usepackage{graphicx}
\aid{0}
\volume{20}
\issue{3}
\pubyear{2014}
\firstpage{1647}
\lastpage{1671}
\doi{10.3150/13-BEJ536} 

\makeatletter
\def\cal{\mathcal}

\def\a{\alpha}
\def\be{\beta}
\def\bX{{\bar X}}
\def\cB{{\cal B}}
\def\const{\mathrm{const.}}
\def\ep{\varepsilon}
\def\etal{\textit{et al.} }
\def\etalc{\textit{et al.} }
\def\ha{{\hat\a}}
\def\hbe{{\hat\be}}
\def\hell{{\hat\ell}}
\def\hj{{\hat\jmath}}
\def\hpi{{\hat\pi}}
\def\hS{{\widehat S}}
\def\htau{{\hat\tau}}
\def\hU{{\widehat U}}
\def\ij{_{ij}}
\def\ka{\kappa}
\def\la{\lambda}
\def\mi{ | }
\def\misrank{_{\mathrm{misrank}}}
\def\AUC{\mathrm{AUC}}
\def\Om{\Omega}
\def\part{\partial}
\def\ra{\to}
\def\rai{\to\infty}
\def\si{\sigma}
\def\tell{{\tilde\ell}}
\def\Th{\Theta}
\def\var{\operatorname{var}}
\def\ze{\zeta}

\makeatother

\begin{document}
\begin{frontmatter}

\title{Feature selection when there are many influential features}
\pdftitle{Feature selection when there are many influential features}
\runtitle{Feature selection}

\begin{aug}
\author[a,b]{\inits{P.}\fnms{Peter}~\snm{Hall}\thanksref{a,b,e1}\ead[label=e1,mark]{halpstat@ms.unimelb.edu.au}},
\author[c]{\inits{J.}\fnms{Jiashun}~\snm{Jin}\thanksref{c,e2}\ead[label=e2,mark]{jiashun@stat.cmu.edu}}
\and
\author[a]{\inits{H.}\fnms{Hugh}~\snm{Miller}\corref{}\thanksref{a,e3}\ead[label=e3,mark]{h.miller@ms.unimelb.edu.au}}
\runauthor{P. Hall, J. Jin and H. Miller} 
\address[a]{Department of Mathematics and Statistics, University of
Melbourne, VIC 3010,
Australia.\\
\printead{e1,e3}}

\address[b]{Department of Statistics, University of California, Davis,
CA 95616, USA.}

\address[c]{Department of Statistics, Carnegie Mellon University,
Pittsburgh, PA 15213,
USA.\\
\printead{e2}}
\end{aug}

\received{\smonth{5} \syear{2011}}
\revised{\smonth{11} \syear{2012}}

%
\begin{abstract}
Recent discussion of the success of feature selection methods has
argued that focusing on a relatively small number of features has been
counterproductive. Instead, it is suggested, the number of significant
features can be in the thousands or tens of thousands, rather than (as
is commonly supposed at present) approximately in the range from five
to fifty. This change, in orders of magnitude, in the number of
influential features, necessitates alterations to the way in which we
choose features and to the manner in which the success of feature
selection is assessed. In this paper, we suggest a general approach
that is suited to cases where the number of relevant features is very
large, and we consider particular versions of the approach in detail.
We propose ways of measuring performance, and we study both theoretical
and numerical properties of the proposed methodology.
\end{abstract}

%
\pdfkeywords{change-point analysis,
classification,
dimension reduction,
feature selection
logit model,
maximum likelihood,
ranking,
thresholding}

\end{frontmatter}

\section{Introduction}\label{sec1}

\subsection{Motivation and summary}\label{sec1.1}
In this paper, we develop statistical methods for determining features
that enable effective discrimination between two populations of very
high dimensional data, when the number of component-wise differences
that provide leverage for discrimination is relatively large but the
sizes of those differences are potentially small. By way of contrast,
conventional approaches to solving this problem tend to rely on
relatively large differences and relatively small numbers of components
where differences occur.

In such problems, it is generally going to be of substantial practical
interest to identify, with reasonable accuracy, the components that
have greatest leverage for correct discrimination. Simply constructing
a classifier, which might depend in a difficult-to-determine way on
differences between two populations, is generally not going to provide
all the information that is sought. However, particularly when the
number of such components is large, we may not be able to identify the
components without error. How accurate can we be, and in what
circumstances is accuracy high? In this paper we shall endeavour to
answer these questions.

Achieving reasonable accuracy can involve relatively computer-intensive
methods, for example, algorithms that need $\mathrm{O}(p^2)$ rather than $\mathrm{O}(p)$
time if the problem is $p$-dimensional. However, if we use an initial,
deterministic dimension reduction step, which decreases dimension to
$q$ where $q\ll p$, then $\mathrm{O}(p^2)$ calculations can be reduced to $\mathrm{O}(p
\log p+q^2)$, where $p \log p$ is the computational cost of ordering
the initial $p$ components. In many cases, we expect $q$ to be a rather
crude upper bound to the true number, $r$ say, of components that
impact on performance of the classifier. The four-stage algorithm that
we shall introduce in Section~\ref{sec2} enables us to reduce computational
expense from $\mathrm{O}(p \log p+q^2)$ to $\mathrm{O}(p \log p+r^2)$. (These
order-of-magnitude calculations ignore the effects of training sample
size, $n$ say, since in the problems we are considering $n$ is
typically much less than $p$, $q$ or $r$ and so has relatively little
impact on the final result.)

Support for the conjecture that $r$ can be quite large, for example, in
genomics problems, has been given by Goldstein \cite{Gol09}, who, in the
words of J.N. Hirschhorn in the same issue of the \textit{New England
Journal of Medicine}, ``builds a speculative mathematical model and
infers that there will be tens of thousands of common variants
influencing each disease and trait'' (Hirschhorn \cite{Hir09}). Goldstein's
\cite{Gol09} calculations are also consistent with $r$ being in the
thousands, not just the tens of thousands:

{\begin{quote}
$\ldots$ the genetic burden of common diseases must be mostly carried
by large numbers of rare variants. In this theory, schizophrenia, say,
would be caused by combinations of 1000 rare genetic variants, not of
10 common genetic variants.
\end{quote}}

\noindent(See Wade \cite{Wad09}.) Kraft and Hunter \cite{KraHun09} argue that ``many, rather
than few, variant risk alleles are responsible for the majority of the
inherited risk of each common disease.'' Again, $r$ is large rather
than small.

In the discussion above one might interpret $p$ and $r$ as representing
numbers of single nucleotide polymorphisms (SNPs), alleles, or perhaps
genes. There are believed to be between 10 and 30 million SNPs on a
human chromosome, and some 25,000 genes. However, genomic analyses
based on decoding the full DNA of individuals who suffer from specific
conditions (Wade~\cite{Wad09}) increase the values of both $p$ and $r$ by
orders of magnitude. (In practice $r$ will be chosen empirically, and
so will actually be a function of the data, but at the level of the
discussion in the present section there is little to be gained by
making this distinction.)

\subsection{Example}\label{sec1.2}
\begin{figure}

\includegraphics{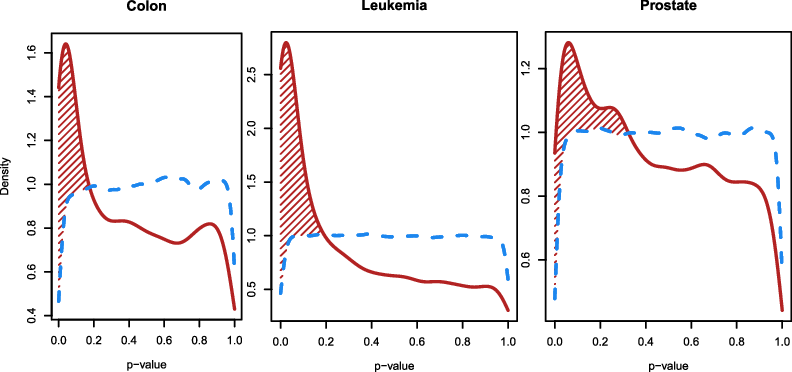}

\caption{Distribution of $p$-values of Wilcoxon tests applied to
the components of microarray datasets. The red solid line represents
the actual data, while the dashed blue line is when the response is
randomly permuted and so there is no genuine signal. The shaded regions
give estimates of the proportion of components genuinely related to the
response.}\label{sim2}
\end{figure}

We do not have to look beyond datasets well covered in the statistical
literature to see evidence of a very large number of influential
components being present in a dataset. We have taken three popular
microarray datasets, described in more detail in Section~\ref{sec5.4},
pertaining to colon cancer, leukemia and prostate cancer. In each case,
there are many thousands of genes, each with a continuum of expression
levels, as well as a binary response variable, indicating whether or
not a patient carries the condition. For each dataset we performed a
Wilcoxon rank sum test for each gene, using the response to divide
observations, so as to test nonparametrically for a change in mean. The
kernel estimates of the densities of $p$-values from these experiments
are displayed in Figure~\ref{sim2} in solid red lines. The peak at the
left for each dataset suggests a significant number of low $p$-values,
more than would be expected if only a few genes contributed to the
understanding of the response.

To extend the experiment further, we randomly permuted the response
amongst the observations 50 times and repeated the Wilcoxon tests,
recovering the distributions denoted by the blue dashed lines. These
can be thought of as the expected distributions of $p$-values were each
gene independently varying with respect to the response. They
asymptotically approximate the uniform distribution. The shaded region
in each plot then provides a rough estimate of the fraction of
components related to the response. These fractions are 16\%, 30\% and
8.5\% of all genes considered (or equivalently, 320, 2140 and 510
genes, resp.), for the colon, leukemia and prostate datasets,
respectively.
We have also done a similar experiment for a SNP data on multiple
sclerosis to be introduced in Section~\ref{sec5.5}, where
the fraction is $3.9\%$ out of 300,900 SNPs.
It should also be noted that these observations are genuine; the
statistic measured (difference between estimated densities) is many
standard deviations above what would be expected by chance.

\subsection{Comments on the literature}\label{sec1.3}
Methods for feature selection based on the linear model are generally
considered only in cases where the number of features is relatively
small. Otherwise, the value of the response variable can be
unreasonably insensitive to changes in a single feature. Examples of
approaches founded on the linear model include the nonnegative garrotte
(e.g., Breiman \cite{Bre95,Bre96}; Gao \cite{Gao98}), the lasso
(Tibshirani \cite{Tib96}), the Dantzig selector (Candes and Tao
\cite{CanTao07}), and related techniques (e.g., Donoho and Huo
\cite{DonHuo01}; Fan and Li \cite{FanLi01,FanLi06}; Donoho and Elad
\cite{DonEla03}; Tropp \cite{Tro05}; Donoho \cite{Don06N1,Don06N2}; Fan
and Ren \cite{FanRen06}; Fan and Fan \cite{FanFan08}). The
feature-ranking approach that we consider is more closely related to
correlation-based approaches of Fan and Lv \cite{FanLv08} and Hall and
Miller \cite{HalMil09}, but it is does not assume the existence of a
response variable. Instead it utilises class labels via a logistic
model. Monograph-length treatments of classifiers and related
methodology include those of Duda \etal\cite{DudHarSto01}, Hastie
\etal\cite{HasTibFri01} and Shakhnarovich \etal\cite{ShaDarInd05}.

\section{Methodology}\label{sec2}

\subsection{Data and algorithm}\label{sec2.1}
Denote the two populations of interest by $\Pi_0$ and $\Pi_1$. Training
data from each are acquired as $p$-vectors $X_i=(X_{i1},\ldots
,X_{ip})$, for $1\leq i\leq n$. We also record the value of a label
$I_i$, for each $i$; it equals 0 or 1, indicating the index of the
population from which $X_i$ came.

One potential algorithm for identifying indices $j$, of vector
components $X\ij$ which capture differences between $\Pi_0$ and $\Pi
_1$, has four stages, (1)--(4) below. The goal of the algorithm is to
determine empirically a set $\{\hj_1,\ldots,\hj_r\}$ of indices, a
subset of $\{1,\ldots,p\}$, such that the features with indices $\hj
_1,\ldots,\hj_r$ have significant influence on whether $X_i$ comes from
$\Pi_0$ or $\Pi_1$. Those features can then be combined into a
classifier, for example, the support vector machine or a centroid-based
method, to effect discrimination.

(1) \textit{Component ranking.}
Using a method such as that suggested in Section~\ref{sec2.2}, rank all
components in terms of their individual influence on $I_i$, interpreted
as a zero--one response variable. This stage takes $\mathrm{O}(p \log p)$ time
to run, and produces a permutation $\hj_1,\ldots,\hj_p$, say, of
$1,\ldots,p$, where the order of the sequence $\hj_1,\ldots,\hj_p$ is
of major importance and signifies that, for each $k$, the component
with rank $\hj_k$ has greater leverage than the component with rank
$\hj
_{k+1}$ on a measure of our ability to predict $I_i$ from $X_i$.

(2) \textit{Deterministic dimension reduction}.
Truncate at $q$ (where $1\leq q\leq p$) the sequence we derived in
step (1). From this point, we work only with $q$-vectors comprised of
the components with indices $\hj_1,\ldots,\hj_q$. The value of $q$ is
determined largely by our computational resources, bearing in mind that
the computational expense of constructing the classifier could be as
high as $\mathrm{O}(q^2)$.

(3) \textit{Adaptive dimension reduction}.
In this stage, we use an empirical method to reduce dimension from $q$,
chosen in stage (3), to $r$, so that the final choice of feature
indices is $\hj_1,\ldots,\hj_r$. Potential approaches are discussed in
Section~\ref{sec2.3}, and include methods based on: (3a) thresholding,
(3b) change-point methods, or (3c) application of classifiers to blocks
of components.

(4) \textit{Backing and filling}.
In practice, it can be advantageous to rerun stage (3) of the algorithm
using several of the values of $j$ chosen early in stage (2), or early
in the implementation of stage~(3), bearing in mind that there is
potential for noise in the choice of $\hj_1$, for example, to throw the
algorithm off course for a period. At this point we could, for example,
experiment with different choices of block size in method (3c).

\subsection{Method for ranking components}\label{sec2.2}
Given an index $j$ between 1 and $p$, and scalar parameters $\a$ and
$\be$, we capture the relationship between $I_i$ and $X\ij$ by assuming
a logit model:
\[
P(I_i=0\mi X\ij) =\bigl\{1+\exp(\a+\be X\ij)\bigr\}^{-1},\qquad
P(I_i=1\mi X\ij)=1-P(I_i=0\mi X\ij) .
\]
The likelihood of $I_i$, given $X\ij$, is
\[
L\ij(I_i\mi\a,\be;X\ij) = \biggl(\frac{t\ij}{1+t\ij}
\biggr)^{ I_i} \biggl(\frac{1}{1+t\ij} \biggr)^{ 1-I_i} =
\frac{t\ij^{I_i}}{1+t\ij} ,
\]
where $t\ij=\exp(\a+\be X\ij)$. Therefore, the negative
log-likelihood is
%
%
\begin{equation}
\ell\ij(\a,\be)=-\log L\ij(I_i\mi\a,\be;X\ij) =-I_i (
\a+\be X\ij) +\log \bigl\{1+\exp(\a+\be X\ij) \bigr\} ,\label{eq2.1}
\end{equation}
and its counterpart for $X_{1j},\ldots,X_{nj}$ is
%
%
\begin{equation}
\ell_j(\a,\be)=\frac{1}{n}\sum_{i=1}^n
\ell\ij(\a,\be) .\label{eq2.2}
\end{equation}
Define $(\ha_j,\hbe_j)$ to be the value of $(\a,\be)$ that minimises
$\ell_j(\a,\be)$, and put
%
%
\begin{equation}
\hell_j=\ell_j(\ha_j,\hbe_j)
.\label{eq2.3}
\end{equation}
The ordering $\hj_1,\ldots,\hj_p$ mentioned in step (1) of the
algorithm in Section~\ref{sec2.1} is determined by the values of $\hell_j$.
Specifically, $\hell_{\hj_1}\leq\cdots\leq\hell_{\hj_p}$.

\subsection{Methods for adaptive dimension reduction}\label{sec2.3}
Several approaches are feasible, including: (3a) \textit{Thresholding.}
Here we compute, from the data, a subsidiary criterion $\tell_j$, for
$1\leq j\leq p$ (we might simply choose $\tell_j\equiv0$) and a
threshold $t$; we take $k_0\geq0$ to be an integer; and we define $r\in
[k_0+1,q]$, a function of the data, to be the least integer in that
range such that $\hell_{\hj_{r+k}}-\tell_{\hj_{r+k}}>t$ for $1\leq
k\leq k_0$. See Section~\ref{sec4} for an example. (3b) \textit{Change-point
methods.} Here we look for a change-point in the sequence $\hell_{\hj
_1},\ldots,\hell_{\hj_p}$, and we take $\hj_r$ to be the location of
that point. (There is a vast literature on methodology and theory for
change-point detection. It includes book-length accounts by Carlstein
\etal\cite{CarMulSie94}, Cs\"org\H o and Horv\'ath \cite{CsoHor97},
Chen and Gupta \cite{CheGup00} and Wu \cite{Wu05}.) (3c)
\textit{Application of classifiers}. For $k\geq1$, let
$\cB_k=\{\hj_{(k-1)b+1},\ldots,\hj_{kb}\}$ denote the $k$th block of
feature indices; here, $b$ denotes block length. (Theoretical
considerations suggest that taking $b\sim\const\ n$ is appropriate.) In
step $s$ of stage (3c) we construct the classifier that is based on the
training data vectors where all but the components with indices in
$\bigcup_{1\leq k\leq s} \cB_k$ have been stripped away. We use
cross-validation to measure classifier performance, and in this way we
determine whether progressing from step $s$ to step $s+1$ gives an
improvement. If it does not then, subject to the ``jiggling'' suggested
in stage (4), we stop at step $s$. If performance is improved by
passing to the $(s+1)$st block then we proceed to step $s+1$, where we
again assess performance. However, this approach can be biased in
favour of low apparent error rate, without having the same impact on
actual error rate; see Section~\ref{sec5.3} for discussion.

\subsection{Duration of algorithm}\label{sec2.4}
Stage (3) of the algorithm takes $\mathrm{O}(r^2)$ time to complete,
where $r$, in the range $1\leq r\leq q$, denotes the final number of
components on which we determine that the classifier should depend. In
particular, the algorithm concludes with a list of $r$ components, say
$\hj_{\ell _1},\ldots,\hj_{\ell_r}$ where
$\{\hj_{\ell_1},\ldots,\hj_{\ell _r}\} =\bigcup_{k\leq{\hat s}}
\cB_k\subseteq\{\hj_1,\ldots,\hj_q\}$, on which the final classifier is
based. The $\mathrm{O}(r^2)$ figure is derived as follows: Constructing
the classifier from $s$ batches takes $\mathrm{O}(sb)$ time, so the
total time needed is $\mathrm{O}(\sum_{1\leq s\leq{\hat s}}
sb)=\mathrm{O}({\hat s}^2b)=\mathrm{O}(r^2)$. Here, since $b$ is
generally of order $n$ (see Section~\ref{sec5}) and ${\hat s}
^2b=\mathrm{O}(r^2b)$, we have replaced $r^2b$ by $r^2$ since, in the
problems we are treating, $n$ is generally so much less than $r$ that
it can be treated as fixed. Provided the number of reruns in stage (4)
is only $\mathrm{O}(1)$, the order of magnitude of the time taken to
run the algorithm to completion is $\mathrm{O}(p \log p+r^2)$.

\subsection{Discussion}\label{sec2.5}
The procedures above are intended to reflect methodologies already used
in practice. Our main aim is to show that such techniques can be used
to address not just contemporary problems where there is believed to be
considerable sparsity and only a small number of significant features
(e.g., five or ten genes out of thousands or tens of thousands), but
also reduced sparsity and a larger total number of features (e.g.,
thousands or tens of thousands of DNA sequences out of tens
or hundreds of thousands of possibilities). Additionally we show that
the methods continue to work well under minimal distributional
assumptions (e.g., normality is not needed), and minimal conditions
about relationships among features. In all these senses, procedures
such as those described above are particularly versatile.

The proposed approach belongs to the class of marginal methods for that
it only uses marginal information for feature selection. Marginal
methods are feasible when different features are independent or weakly
dependent. When this is not the case,
a natural problem is how to incorporate dependence among different
features for feature ranking and classification. One possible approach
is to first estimate the correlations and then include them in the
classifier. However, without rather narrow assumptions this can
introduce substantial noise into the process of inference for
high-dimensional data, and in fact one can be better off simply
ignoring the correlations; see, for example, Tibshirani \etal\cite{Tibetal02},
Bickel and Levina \cite{BicLev04} and Fan and Fan \cite{FanFan08}.

The feature ranking problem can also be re-cast in terms of logistic
regression, where the design matrix $\{X_{ij}\}_{\{1 \leq i \leq n, 1
\leq j \leq p\}}$ is randomly generated and perhaps non-Gaussian.
Hence, there is a natural connection between the
component-by-component, or marginal, approach suggested in the present
paper, and existing techniques for feature selection, such as the lasso
(e.g., Chen \etalc\cite{CH98}; Tibshirani \cite{Tib96}). Genovese \etal\cite{Genetal12}
compared marginal methods with the lasso and showed that the two
approaches have comparable performance over a wide range of parameter
values, both theoretically and empirically, but that marginal methods
are computationally more efficient, especially when both $n$ and $p$
are large.

Our work is connected to that of Fan and Lv \cite{FanLv08} in that it is
founded on component-wise analysis of high-dimensional data. Fan and
Lv's ingenious technique was based on correlation, and is related to
likelihood through its connections to least-squares in the setting of
Gaussian experimental errors. In contrast, we use a logistic model to
define a likelihood for each component, and we rank those
component-wise likelihoods. Fan and Lv's \cite{FanLv08} aim is ``to reduce
dimensionality from high to a moderate scale that is below the sample
size'' (to quote from their paper), and then to use relatively
conventional methods to select features. However, our goal is
inherently different, since a basis for our argument is the concern
that, as Goldstein \cite{Gol09} argued, the number of important features
might be relatively large, and in particular much larger than sample size.

Our paper is connected to those of Fan and Lv \cite{FanLv08} and Fan and Song
\cite{FanSon10}, and in particular our methodology, minus the ``adaptive
dimension reduction'' and ``backing and filling'' steps, is Fan and
Song's \cite{FanSon10} SIS approach when the model in question is logistic.
However, the parameter settings that we treat, both in our numerical
work and our theory, are quite different from theirs. In particular, we
focus on regimes where the useful features are very weak and only
moderately sparse. In contrast, Fan and Lv \cite{FanLv08} and Fan and Song
\cite{FanSon10} address cases where the signals are relatively strong or very
sparse. Secondly, we evaluate our approach using the new criterion of
misranking, introduced below, which is particularly appropriate to
problems where many influential features are present and, for at least
this reason, different from criteria employed by Fan and Lv \cite{FanLv08} and
Fan and Song \cite{FanSon10}. It is also different from the FDR approach adopted
by Benjamini and Hochberg \cite{BenHoc95} and many other authors.

In broad terms the methodology that we suggest could be used in
conjunction with existing approaches to feature selection, but those
generally do not conclude with a ranking of features, and that would
limit their usefulness to us. More particularly, as discussed above,
our approach differs from those of Fan and Lv \cite{FanLv08} and others, in
that our aim is not to dramatically reduce the number of features to a
handful, but to address a much larger number of features. Nevertheless
some methods, for example, random forests (which produce a feature
importance list) and shrunken centroids do create an ordering for the
portion of predictors that are given nonzero weight, and would be
candidates for implementing part of our methodology.

\section{Properties of feature ranking}\label{sec3}

\subsection{Main result on ranking}\label{sec3.1}
Let $\pi=E(I_i)$ denote the proportion of data that come from
population $\Pi_1$. We assume below that $0<\pi<1$, and that when the
training data are drawn randomly from the union of $\Pi_0$ and $\Pi_1$,
the prior probability that any given datum $X_i$ is from $\Pi_1$ equals
$\pi$. Therefore, the corresponding probability for $\Pi_0$ is $1-\pi$.
We take $n$, representing the total size of the training sample, to be
the key asymptotic parameter, and interpret the dimension, $p$, of
$X_i$ as a function of $n$.

Next, we describe our model. Write $X_i=(X_{i1},\ldots,X_{ip})$, and,
when $X_i$ comes from $\Pi_1$, take
%
%
\begin{equation}
X\ij=\mu_j+Z\ij,\label{eq3.0}
\end{equation}
where $\mu_j\geq0$ are constants and, for simplicity, we assume that
the variables $Z\ij$ are identically distributed as $Z$, say. (This
condition can be relaxed; see (\ref{eq3.4}) below.) Take $X\ij=Z\ij$ when $X_i$
is drawn from $\Pi_0$. The vectors $(Z_{i1},\ldots,Z_{ip})$ and
variables $I_i$, for $1\leq i\leq n$, are assumed to be totally
independent. We allow the $\mu_j$'s to be functions of $n$.

Write $\ha_j$ and $\hbe_j$ for the values of $\a$ and $\be$ that
jointly minimise $\hell_j(\a,\be)$ within radius $n^{-c}$ of $(\a
_0,0)$, where $\a_0=\log\{\pi/(1-\pi)\}$, for any given $c\in
(0,{ {\frac{1}{2}}})$.
(A separate argument can be used to prove that if $\ha_j$ and $\hbe_j$
are chosen without constraint then, under the conditions of Theorem~\ref{th1}
below, they satisfy $\ha_j=\a_0+\mathrm{O}_p(n^{-c})$ and $\hbe_j=\mathrm{O}_p(n^{-c})$
uniformly in $1\leq j\leq p$, for some $c\in(0,{ {\frac{1}{2}}})$.) Define
%
%
\begin{equation}
S_j=\frac{\sum_i (I_i-\pi) X\ij}{ n \pi(1-\pi)} ,\label{eq3.1}
\end{equation}
and note that $E(S_j)=0$. Put $\la=(n^{-1}\log n)^{1/2}$.

\begin{thm}\label{th1} Assume that for each $n$, $|\mu_j|\le\const\ \la$ for
$1\leq
j\leq p$; that $p=p(n)\rai$ and, for constants $B_1>0$ and $B_2>2\max
(B_1+3,2B_1)$, $p=\mathrm{O}(n^{B_1})$ and $0<E|Z|^{B_2}<\infty$; and that
$E(Z)=0$. Then, uniformly in $1\leq j\leq p$,
%
%
\begin{equation}
\hell_j=\ell_j(\ha_j,\hbe_j)
=R-{ {\tfrac{1}{2}}}\pi(3-2 \pi) \bigl(EZ^2
\bigr)^{-1}(S_j+\mu_j)^2
+\mathrm{O}_p \bigl(\la^3 \bigr) ,\label{eq3.2}
\end{equation}
where the random variable $R=R(n)$ does not depend on $j$. More
particularly, the $\mathrm{O}_p(\la^3)$ term in (\ref{eq3.2}) can be written as
$\Th_j \la^3$, where, for a constant $B_3$ depending on $B_1$ and
$B_2$, and with $B_4={ {\frac{1}{4}}}\{B_2-2\max(B_1+3,2B_1)\}-\ep$
for any $\ep
>0$, the random variables $\Th_j$, for $1\leq j\leq p$, satisfy:
%
%
\begin{equation}
\sum_{j=1}^p P\bigl(|\Th_j|>B_3\bigr)=\mathrm{O}
\bigl(n^{-B_4} \bigr)\label{eq3.3}
\end{equation}
as $n\rai$.
\end{thm}

The statistic $S_j$ is, up to normalisation, the well-known $z$-score
statistic for testing whether the $j$th feature is significant or not;
see, for example, Donoho and Jin \cite{DonJin08,DonJin09} and Jin \cite{Jin09}. In~(\ref
{eq3.2}), since the first term, $R$, does not depend on $j$ then the
second term is the one that reflects the strengths of individual
features. As a result, ranking features according to $\hell_j$ gives
close, but not necessarily the same, results as ranking features
according to $|S_j+\mu_j|$.

Next, we interpret the assumptions imposed in Theorem \ref{th1}. The condition
$|\mu_j|\le\const\ \la$, for $1\leq j\leq p$, is imposed so as to make
the contribution of $\mu_j$ difficult to identify. In particular, if
$|\mu_j|$ is of larger order than $\la$ then it can be shown from large
deviation properties that the contribution of the ``signal'' $\mu_j$
will be easily visible above the noise, and so the contribution of the
gene corresponding to the index $j$ will be very easy to identify. This
would exemplify the classical setting, where a relatively small number
of genes describe adequately the link to disease, and we are not
addressing that case in this paper. The conditions $p=\mathrm{O}(n^{B_1})$ and
$0<E|Z|^{B_2}<\infty$ imply that the number of dimensions can be
polynomially large as a function of sample size, and that only a
polynomial number of moments needs to be assumed. The assumption
$B_2>2\max(B_1+3,2B_1)$ implies that the number of moments required
increases no faster than linearly as a function of the rate of growth
of $p$, expressed in terms of sample size.

The assumption that each $Z\ij$ has the same distribution is made to
simplify discussion and notation, and is readily relaxed. For example,
it suffices to assume that each $Z\ij$, for $1\leq i\le n$, is
distributed as $Z_j$, say, where, instead of the assumptions imposed on
$Z$ in the theorem, we ask that:
%
%
\begin{equation}\label{eq3.4}
\begin{tabular}{p{280pt}} $E(Z_j)=0$, $E
\bigl(Z_j^2\bigr)$ is bounded away from zero and infinity
uniformly in $j$, and $P(|Z_j|>x)\leq P(|Z|>x)$ for all $x>0$ and
some random variable $Z$, where $E|Z|^{B_2}<\infty$.
\end{tabular}
\end{equation}
In this case, the moment $E(Z^2)$ in (\ref{eq3.2}) would be replaced
by $E(Z_j^2)$. The conclusions that we draw, below, from Theorem \ref{th1} are
unchanged, provided we interpret $\mu_j$ as $\mu_j/(EZ_j^2)^{1/2}$
during discussion.

Although we ask that the vectors $X_i$ be independent, we make no
assumption about the relationships among their components. For example,
the values of $Z_{i1},\ldots,Z_{ip}$ can be highly dependent (indeed,
in an extreme case, equal to one another) or completely independent.
The latter instance is actually the most difficult, in terms of
rigorously establishing that (\ref{eq3.2}) and (\ref{eq3.3}) hold. At
the other end of the spectrum, the case where $Z_{i1}=\cdots=Z_{ip}$
with probability 1 is trivial, since there effectively only a single
component index, with different candidate values for the mean, has to
be treated.

\subsection{Expected number of misrankings}\label{sec3.2}
Assume that some of the $\mu_j$'s are zero and all the others are
strictly positive. Ideally, we would like the criterion $\hell_j$ to be
a good indicator of the positivity of $\mu_j$, in particular to take a
lesser (or larger negative) value if $\mu_j$ is positive than it does
when $\mu_j=0$. Reflecting this aspiration, if there exist component
indices $j_1$ and $j_2$ such that $\mu_{j_1}>0$ and $\mu_{j_2}=0$, but
$\hell_{j_2}<\hell_{j_1}$, then we shall say that a misranking has
occurred. The expected total number of misrankings,
\[
\nu\misrank=\sum_{j_1 \dvt \mu_{j_1}>0} \sum
_{j_2 \dvt \mu_{j_2}=0} P(\hell_{j_2}<\hell_{j_1}) ,
\]
is a measure of the performance of $\hell_j$ as a criterion for
distinguishing between positive and zero values of $\mu_j$; lower
values of $\nu\misrank$ correspond to higher performance.

Since the random variable $S_j$ in (\ref{eq3.1}) and (\ref{eq3.2}) has
standard deviation of size $n^{-1/2}$ then, if the positive $\mu_j$'s are
of smaller order than $n^{-1/2}$, with probability converging to ${
{\frac{1}{2}}}$
any attempt to rank any pair of means $\mu_j$ using the values of
$\hell
_j$ will produce\vadjust{\goodbreak} the wrong result about half the time. The following
theorem makes this clear. Let $p_1$ denote the number of indices $j$
for which $\mu_j>0$, and put
%
%
\begin{equation}
\si_{j_1j_2,\pm}^2 =\pi(1-\pi) \var(Z_{1j_1}\pm
Z_{1j_2})\label{eq3.5}
\end{equation}
for respective choices of the plus and minus signs.

To avoid degeneracy, where the problem of identifying genes is
relatively simple, we assume below that the quantities $\si
_{j_1j_2,\pm
}^2$ are bounded away from zero, and of course, since we are
endeavouring to capture cases where the signals $\mu_j$ are so weak
that the corresponding genes cannot be identified, we ask that the $\mu
_j$'s be uniformly of smaller order than $n^{-1/2}$.

\begin{thm}\label{th2} Assume the conditions of Theorem \ref{th1}, that $\sup_{j\leq p}
\mu_j=\mathrm{o}(n^{-1/2})$, and that $\si_{j_1,j_2,\pm}^2$ is bounded away from
zero uniformly in $j_1$, $j_2$ and choices of the $\pm$ signs. Then
$P(\hell_{j_2}<\hell_{j_1})\ra{ {\frac{1}{2}}}$ uniformly in $1\leq
j_1,j_2\leq p$,
in particular uniformly in pairs $j_1,j_2$ such that $\mu_{j_1}>0$ and
$\mu_{j_2}=0$. Moreover, $\nu\misrank={ {\frac{1}{2}}}p_1
(p-p_1)+\mathrm{o}\{p_1
(p-p_1)\}$ as $n\rai$.
\end{thm}

Likewise, if the positive $\mu_j$'s are of size $n^{-1/2}$ then the
probability of incorrectly ranking the $j_1$th component lower than the
$j_2$th component, even though $\mu_{j_1}>0$ and $\mu_{j_2}=0$, does
not converge to zero. The next theorem quantifies this property. There
we define $\Phi$ to be the standard normal distribution function.

\begin{thm}\label{th3} Assume the conditions of Theorem \ref{th1}, and that each nonzero
$\mu_j$ equals $c n^{-1/2}$ where $c>0$. Then
\[
P (\hell_{j_2}<\hell_{j_1} ) =\Phi(-c/\si_{j_1j_2,+})
\Phi(c/\si_{j_1j_2,-}) +\Phi(c/\si_{j_1j_2,+}) \Phi(-c/
\si_{j_1j_2,-})+\mathrm{o}(1) ,
\]
uniformly in $j_1,j_2$ such that $\mu_{j_1}>0$ and $\mu_{j_2}=0$. Furthermore,
\begin{eqnarray*}
\nu\misrank&=& \sum_{j_1 \dvt \mu_{j_1}>0} \sum
_{j_2 \dvt \mu
_{j_2}=0} \bigl\{\Phi(-c/\si_{j_1j_2,+}) \Phi(c/
\si_{j_1j_2,-})
 +\Phi(c/\si_{j_1j_2,+}) \Phi(-c/\si_{j_1j_2,-}) \bigr\} \\
&&{}
+\mathrm{o}\bigl
\{p_1 (p-p_1)\bigr\}
\end{eqnarray*}
as $n\rai$.
\end{thm}

Here, by default, $\Phi(-\infty)=0$ and $\Phi(\infty)=1$. This is
relevant when $\si_{j_1 j_2,\pm} =0$.

If the number of components where the mean is positive is large, for
example if it equals a nonnegligible proportion of the total number,
$p$, of components, then the number of misrankings can generally not be
reduced to low levels unless we take the nonzero means to be a little
larger than $n^{-1/2}$ in order of magnitude terms. It is enough to take
the positive mean to be a logarithmic factor larger; specifically, the
mean should equal $c \la$ where, as before, $c>0$ and $\la
=(n^{-1}\log
n)^{1/2}$. Theorem \ref{th4}, below, shows that in this case the expected number
of misrankings can be reduced to a quantity of smaller order than $p$,
or even to a number that converges to zero polynomially fast, depending
on how large we choose $c$.

As a prelude to stating Theorem \ref{th4}, we introduce an assumption which
asks that the random variables $Z_j$ satisfy a pairwise Cram\'er
continuity condition. Here it is convenient to assume that there exists
an infinite stochastic process $Z_1,Z_2,\ldots$ such that:

%
\begin{equation}\label{eq3.6}
\begin{tabular}{p{280pt}} (a) Each $Z_j$ has the
distribution of $Z$ (in this sense the process $Z_1,Z_2,
\ldots$ is weakly stationary), (b) if $X_i=(X_{i1},
\ldots,X_{ip})$ is drawn from $\Pi_0$ then
$X_{i1},\ldots,X_{ip}$ has the same joint distribution as
$Z_1,\ldots ,Z_p$, and (c) if $X_i$ is drawn
from $\Pi_1$ then $X_{i1},\ldots,X_{ip}$ has the
same joint distribution as $Z_1+\mu_1,
\ldots,Z_p+\mu_p$, where the $\mu_j$'s are the
nonnegative constants introduced prior to Theorem \ref{th1}.
\end{tabular}
\end{equation}

The Cram\'er continuity condition we impose is the following:
%
%
\begin{equation}
\limsup_{t\rai} \sup_{|t_1|+|t_2|>t} \sup
_{1\leq j_1<j_2<\infty} \bigl|E\bigl\{\exp(\mathrm{i}t_1 Z_{j_1}+\mathrm{i}t_2
Z_{j_2})\bigr\}\bigr| <1 ,\label{eq3.7}
\end{equation}
where on this occasion $\mathrm{i}=\sqrt{-1}$. For example, (\ref{eq3.7}) would
hold if the process $Z_j$ were strictly stationary and each pair
$(Z_{j_1},Z_{j_2})$ had a joint density $f_{j_1j_2}$ that satisfied
$\sup_{j_1,j_2}\int\int|{\ddot f}_{j_1j_2}|<\infty$, where ${\ddot
f}_{j_1j_2}(x_1,x_2)=(\part^2/\part x_1\part x_2)
f_{j_1j_2}(x_1,x_2)$. It would also hold if the variables $Z_j$ were
independent with a common nonsingular distribution.

Recall that $p_1$ equals the number of indices $j$ such that $\mu_j>0$,
and that $B_4=B_4(\ep)={ {\frac{1}{4}}}\{B_2-2\max(B_1+3,2B_1)\}-\ep
$ where
$\ep
>0$. Define $\si_{j_1j_2,\pm}$ by (\ref{eq3.5}) and put $\ka
_n=(\log
n)^{1/2}$.

\begin{thm}\label{th4} Assume the conditions of Theorem \ref{th1}, that $(\ref{eq3.6})$
and $(\ref{eq3.7})$ hold, and that $B_2$, in the moment condition
$E|Z|^{B_2}<\infty$, is so large that for some $\ep>0$, $p_1
(p-p_1)=\mathrm{o}(n^{B_4})$. Take each nonzero $\mu_j$ to equal $c (n^{-1}\log
n)^{1/2}$, where $c>0$. Then

%
\begin{eqnarray}
\label{eq3.8} \nu\misrank &=&\bigl\{1+\mathrm{o}(1)\bigr\} \sum
_{j_1 \dvt \mu_{j_1}>0} \sum_{j_2 \dvt \mu_{j_2}=0} \bigl\{\Phi(-c
\ka_n/\si_{j_1j_2,+})
+\Phi(-c \ka_n/\si_{j_1j_2,-})\bigr\}
\nonumber
\\[-8pt]
\\[-8pt]
\nonumber
&&{} +\mathrm{o}(1)
\end{eqnarray}
as $n\rai$.
\end{thm}

Elucidation of (\ref{eq3.8}) requires information about the covariance
of the process $Z_j$, in (\ref{eq3.6}). For simplicity let us assume
that the variables $Z_j$ are uncorrelated. Then by (\ref{eq3.5}), $\si
_{j_1j_2,\pm}^2=2 \pi(1-\pi) E(Z^2)\equiv(2c_0)^{-1}$, say, for
either choice of the $\pm$ signs. Hence, (\ref{eq3.8}) implies that
\[
\nu\misrank =\bigl\{1+\mathrm{o}(1)\bigr\} p_1 (p-p_1) \bigl\{2
\pi c (2 c_0 \log n)^{1/2} \bigr\}^{-1}n^{-c_0c^2}
+\mathrm{o}(1) .
\]
Therefore, if $c$ is chosen so large that $p_1 (p-p_1) (\log n)^{-1/2}
n^{-c_0c^2}\ra0$ then the expected number of misranks will converge to
zero. For smaller positive values of $c$ the expected number will be of
smaller order than the potential number of misranks, $p_1 (p-p_1)$,
but it will not necessarily be negligible itself.

The results in Theorems \ref{th2}--\ref{th4} have benefited from a simplification
afforded by the assumption that the variables $Z\ij$ all have the same
distribution, and in particular have the same variance. As noted below
Theorem \ref{th1}, that condition can be relaxed and the assumption (\ref
{eq3.4}) imposed instead. In practice, however, one could standardise,
in a componentwise fashion, the values of $X\ij$ for scale, and in that
case it is possible to state versions of Theorems \ref{th2}--\ref{th4} in settings
where $E(Z\ij^2)$ varies with $j$. The model that we have been using,
that is, $X\ij=Z\ij+c (n^{-1}\log n)^{1/2}I_i$ where the $Z\ij$'s are
independent, is (for moderate $n$) a good approximation to the
standardised form $X\ij'=Z\ij'+c (n^{-1}\log n)^{1/2}I_i$, where the
$Z\ij'$s satisfy (\ref{eq3.4}). Detailed arguments here are similar to
those given by Hall and Wang \cite{HalWan10}.

\subsection{Effects of dependence of the process $Z_j$ on
interpretations of (\texorpdfstring{\protect\ref{eq3.2}}{3.3})}\label{sec3.3} The expected value of
the number of misrankings, which we treated in Section~\ref{sec3.3}, is not as
much affected by dependence among components of the $Z_j$ process as
are other aspects of the distribution of the number of misrankings. For
example, if the $Z_j$'s (in the stochastic process $Z_1,Z_2,\ldots$
introduced in (\ref{eq3.6})) are all independent then the quantities
$S_j$, defined at (\ref{eq3.1}) and on which the values of $\hell_j$
predominantly depend (see (\ref{eq3.2})), are also independent, and so
decisions based on the respective values of $\hell_j$ are made
virtually independently of one another. In this case the variance of
the total number of misrankings is\vadjust{\goodbreak} relatively low. However, if the
$Z_j$'s are highly dependent then the variance can be higher, although
it depends on how the positive means are distributed among the
components of $X_i$. In the present section we briefly discuss these
issues.

Let the process $Z_1,Z_2,\ldots$ in (\ref{eq3.6}) be $\ze$-dependent,
meaning that any subsequence\break  $Z_{j_1},\ldots,Z_{j_k}$ such that
$j_{\ell
+1}-j_\ell>\ze$, for each $\ell$, is comprised entirely of independent
random variables. We permit $\ze$ to diverge with $n$, and we suppose
that $p_1$, the number of nonzero values of $\mu_j$, can also increase
with $n$ and that $\limsup_{n\rai} p_1/p<1$. One approach to arranging
the nonzero means is to distribute them randomly, for example, taking
$\mu_j=J_j \mu$ for $1\leq j\leq p$, where $J_1,\ldots,J_p$ is a
random permutation of $p_1$ ones and $p-p_1$ zeros and is independent
of the $Z_j$'s in (\ref{eq3.6}). In this setting, the clustering that
arises through dependence is often negligible, even if the dependence
in the process $Z_1,\ldots,Z_p$ is quite strong.

To appreciate why, note that the expected number of nonzero means in
each string of $\ze$ consecutive components of $X_i$, when $X_i$ is
drawn from $\Pi_1$, equals $\ze\times p^{-1}p_1$; and that if $(\ze
p_1)^2=\mathrm{o}(p)$ then the probability that none of the approximately $p/\ze
$ strings (placed end to end) of $\ze$ consecutive components that
contain one or more nonzero means are adjacent, and the probability
that none of the strings contains more than one component, both
converge to zero. Therefore, in view of the assumption of $\ze
$-dependence, if we treat the $Z_j$'s as independent and identically
distributed when making a statement about properties of rankings
deduced from (\ref{eq3.2}), the probability that we commit an error in
the statement converges to zero as $n\rai$. It can then be deduced
that, in cases where the positive means are randomly distributed and
$(\ze p_1)^2=\mathrm{o}(p)$, the variance of the number of misrankings is
relatively low.

Alternatively, rather than scatter the nonzero means $\mu_j$ randomly
throughout the vector $(Z_1,\ldots,Z_p)$, we could place them all down
one end. This makes the distribution of those quantities just about as
``clumpy'' as possible, by exploiting the $\ze$-dependence property.
For example, if $p_1\leq\ze$ then all of the nonzero means are
attributed to the first $p_1$ variables in the sequence $X_{i1},\ldots
,X_{ip}$, when $X_i$ is drawn from population $\Pi_1$. The assumption
of $\ze$-dependence permits $X_{i1}=\cdots=X_{ip_1}$ with
probability 1, whenever $X_i$ comes from either $\Pi_0$ or $\Pi_1$.
This reduces the amount of available information, since all the
components that contain information for discriminating between $\Pi_0$
and $\Pi_1$ are simply copies of one another; there are no independent
sources of corroborating information. Moreover, the values of $S_j$ in
(\ref{eq3.2}) are identical for $1\leq j\leq p_1$, and so the values of
$\hell_j$ are the same too, up to remainders of order $\la^3$, which
implies that the total number of misranks is approximately equal to
$p_1$ times the number of times that a specific component with a
positive mean is misranked.

Of course, this increases the variance of the number of misranks. The
setting $p_1\leq\ze$ can encompass instances where $(\ze p_1)^2=\mathrm{o}(p)$,
which was shown two paragraphs above to result in a relatively high
amount of information about the differences between $\Pi_0$ and $\Pi_1$
when the nonzero means are scattered randomly in the data vector. These
examples illustrate the more general rule that, in cases where the
positive means are distributed consecutively in relatively long-range
dependent vectors $X_i$, the variance of the number of misranks tends
to be higher than in cases where those means are distributed at random.

We conclude this section by comparing the notion of misranking with
that of (feature) False Discovery Rate (FDR); for the latter, see for
example Benjamini and Hochberg \cite{BenHoc95} and\vadjust{\goodbreak} Abramovich \textit{et al.} \cite{Abretal06}. For
any set of selected features, FDR equals the fraction of falsely
selected features. This concept and that of misranking both provide
informative measures of how well important features are ranked, but
they are nevertheless different in important respects. To appreciate
why, let us focus on a specific feature. If we adhere to the notion of
FDR then all that matters is whether the feature is selected or not. If
we instead we use the concept of misranking then the order or rank of
the feature being selected also matters. Technically there are also
important differences. For instance, misrankings are defined quite
simply in terms of pairwise comparisons of individual features, while
FDR can involve higher-order relationships among different features. If
we consider the influence, on these measures, of dependence among
features, then misranking depends only on pairwise dependence, but FDR
may depend on high-order relationships. As a result, FDR can be
significantly more difficult to characterise than misranking, and
requires much more heavily constrained assumptions about dependence
than are necessary using the misranking measure. Therefore, since the
central problem is how well important features are ranked, it is more
appropriate to assess performance here using misranking, rather than FDR.

The number of misrankings bears a close relationship to both the
Wilcoxon rank-sum test and the area under curve ($\AUC$) of the ROC plot
(see Hanley and McNeil \cite{HanMcn82}). In this case, the ROC is constructed
with respect to whether each feature is correctly classified as having
nonzero mean on not. In fact, it is possible to show that
%
%
\begin{equation}
\AUC= 1 - \bigl\{p_1(p-p_1)\bigr\}^{-1} (\rm{\#
misrankings}) .\label{eq3.9}
\end{equation}
Thus, we may interpret properties of $\nu\misrank$ in Theorems \ref{th2}--\ref{th4} in
the context of $\AUC$. For instance, under the assumptions of Theorem \ref{th2}
the $\AUC$ score decays to 0.5, this being the score for the random
guessing model.

\section{Thresholding for adaptive dimension reduction}\label{sec4}

Recall that in Theorem \ref{th1} we showed that $\hell_j$ equals
$-U_{j1}$, where
%
%
\begin{equation}
U_{j1}={ {\tfrac{1}{2}}}\pi(3-2 \pi) \bigl(EZ^2
\bigr)^{-1}(S_j+\mu _j)^2 ,
\label{eq4.1}
\end{equation}
plus a quantity that does not depend on $j$, plus a remainder term that
is uniformly smaller in size. The feature-ranking step (stage (1) in
the algorithm in Section~\ref{sec2.1}) aspires to re-order the indices~$j$ such
that the indices for which $\mu_j\neq0$ are ranked first, and those for
which $\mu_j=0$ are listed together at the end of the sequence. If this
objective is largely achieved, then the main task that remains is to
choose the point in the ranking where the change occurs; this is
stage (3) of the algorithm in Section~\ref{sec2.1}. In the present section we
explore method (3a), based on thresholding; see Section~\ref{sec2.3}.

Observe that if $U_{j1}$ is as at (\ref{eq4.1}) then
$U_{j1}=U_{j2}+U_{j3}$, where
%
%
\begin{equation}
U_{j2}={ {\tfrac{1}{2}}}\pi(3-2 \pi) \bigl(EZ^2
\bigr)^{-1}S_j^2 ,\qquad U_{j3}={ {
\tfrac{1}{2}}}\pi(3-2 \pi) \bigl(EZ^2 \bigr)^{-1} \bigl(
\mu _j^2+2 \mu_j S_j\bigr) .
\label{eq4.2}
\end{equation}
If we can construct a good approximation, $\hU_{j2}$ say, to $U_{j2}$
then we can subtract it from $\hell_j$, leaving only $U_{j3}$ plus
a\vadjust{\goodbreak}
small remainder. The value of $U_{j3}$ is exactly zero if $\mu_j=0$,
and is strictly positive with high probability if $\mu_j>0$. The
quantity $\tell_j=-\hU_{j2}$ is referred to in that notation in
method (3a) in Section~\ref{sec2.3}. If we choose an appropriate threshold, $t$
say, then we can implement (3a) as follows:
%
%
\begin{equation}\label{eq4.3}
\begin{tabular}{p{280pt}} define $r\in[k_0+1,q]$, a
random variable, to be the least integer in that range such that $
\hell_{\hj_{r+k}}-\tell_{\hj_{r+k}}>t$ for $1\leq k\leq k_0$,
\end{tabular}
\end{equation}
where $k_0\geq0$ is a fixed integer. Then, subject to the jiggling step
in stage (4) of the algorithm, we determine that the features with
indices $\hj_1,\ldots,\hj_r$ are the ones that have greatest influence
on whether a data value $X_i$ came from $\Pi_0$ or $\Pi_1$.

To define $\hU_{j2}$, put $\hpi=n^{-1}\sum_i I_i$ and $\bX_j=n^{-1}
\sum_i
X\ij$, define our estimator\vspace*{1pt} of $\tau^2=E(Z^2)$ by $\htau^2=(np)^{-1}
\sum_i \sum_j (X\ij-\bX_j)^2$, and let $\hS_j=\{n \hpi(1-\hpi)\}^{-1}
\sum_i (I_i-\hpi) (X\ij-\bX_j)$; compare the definition of $S_j$
at (\ref
{eq3.1}). Then, motivated by the definition of $U_{j2}$ at (\ref
{eq4.2}), put
\[
-\tell_j=\hU_{j2}\equiv{ {\tfrac{1}{2}}}\hpi(3-2
\hpi) \htau ^{-1/2}\hS_j^2 .
\]
Theorem \ref{th5}, below, shows in effect that this is a good approximation
to $U_{j2}$. Let $B_4$ be as in Theorem \ref{th1}.

\begin{thm}\label{th5} Under the conditions of Theorem \ref{th1}, we can write
%
%
\begin{equation}
\hell_j-\tell_j=R-U_{j3}+\Om_j
\la^3 ,\label{eq4.4}
\end{equation}
where $R$ is as in $(\ref{eq3.2})$ and, for a constant $B>0$, the
random variables $\Om_j$ satisfy
%
%
\begin{equation}
\sum_{j=1}^p P\bigl(|\Om_j|>B\bigr)=\mathrm{O}
\bigl(n^{-B_4} \bigr) .\label{eq4.5}
\end{equation}
\end{thm}

Finally, we describe implementation of method (3a) in stage (3) in
Section~\ref{sec2.3}. We shall assume we are in the context of Theorem \ref{th4}, where
the nonzero $\mu_j$'s equal $c (n^{-1}\log n)^{1/2}$. Here it is
appropriate to take the threshold $t=t(n)$ to be a sequence of negative
numbers such that $|t|$ is of strictly larger order than $\la^3$ and of
strictly smaller order than $\la^2$; that is, such that
%
%
\begin{equation}
|t|=\mathrm{o} \bigl(\la^2 \bigr)\quad \mbox{and}\quad \la^3=\mathrm{o}\bigl(|t|\bigr) .
\label{eq4.6}
\end{equation}

Let $k_1$ denote the number of nonzero means added to the components of
$X_i$ when $X_i$ is drawn from $\Pi_1$. To simplify discussion, we
assume that $k_1$ is nonrandom, although of course it may depend
on $n$. Below Theorem \ref{th4} we discussed a case where the expected number
of misrankings converged to zero, and hence the probability of a
misranking occurring also tended to zero. In this setting,

%
\begin{equation}\label{eq4.7}
\begin{tabular}{p{280pt}} with probability converging to 1, $
\mu_{\hj_k}>0$ for $k\leq k_1$ and $\mu_{\hj_k}=0$ for
$k>k_1$,
\end{tabular}
\end{equation}
and it can be shown that if $t<0$ satisfies (\ref{eq4.6}) then the
definition of $r$ at (\ref{eq4.3}) produces a random variable which,
with probability converging to 1 as $n\rai$, equals $k_0+k_1$.
Therefore, taking $k_0=0$ in the rule at (\ref{eq4.3}) ensures that,
with probability converging to 1, $r$ is exactly equal to~$k_1$. Cases
where the nonzero means are of size $c (n^{-1}\log n)^{1/2}$, but $c$
is not sufficiently large to ensure that (\ref{eq4.7}) holds, can be
treated satisfactorily by choosing $t<0$ to satisfy (\ref{eq4.6}) but
taking $k_0\geq1$. Depending on the strength of dependence between
components of the vector $X$ it is possible to choose a fixed $k_0$
such that, with probability converging to 1, $r$ is within $C\rho k_1$
of $k_1$, where $C>0$ and $\rho$ equals the expected proportion of
feature indices that have positive means when $X_i$ is drawn from $\Pi
_1$, but are incorrectly ranked at a low level.


\section{Numerical properties}\label{sec5}
Our numerical studies involve three parts. First, in Sections~\ref{sec5.1}--\ref{sec5.2},
we use
simulated data to investigate several aspects of misranking,
including misranking stability, influence of dependence on misranking, and
prediction with large simulated examples.
Second, in Section~\ref{sec5.4}, we compare our methods with several popular
classifiers (including SVM and Random Forest) with the three
aforementioned gene microarray data sets.
Last, in Section~\ref{sec5.5}, we further compare our method with SVM and Random
Forest with
a multiple sclerosis SNP data set.

\subsection{Stability of misranking totals}\label{sec5.1} Here we
simulate under the model at (\ref{eq3.0}), where the signals are
represented by $\mu _j$ and the noise by $Z\ij$. If we measure
performance in terms of the total number of misrankings, or
equivalently in terms of $\AUC$ (see (\ref{eq3.9})), and if the
$\mu_j$'s decrease at rate $n^{-1/2}$, then Theorem \ref{th3} implies that
performance should be stable as a function of sample size, $n$. That
is, it should depend very little on $n$. To explore this property
numerically we consider the cases $n=20$, 50, 100 and 200, with $p= 0.4
\times n^2$.
We take 90\% of the $\mu_j$'s to equal zero and the others to equal
$c_0 (20/n)^{1/2}$, where $c_0=0.4, 0.8, \ldots, 0.16$. The noise
variables $Z\ij$ are independent and identically distributed as N$(0,1)$.

\begin{figure}

\includegraphics{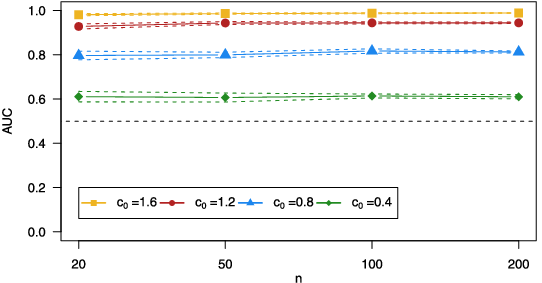}

\caption{AUC scores in simulation study with $\mu_j = c_0
(20/n)^{1/2}$, in the example of Section~\protect\ref{sec5.1}.}\label{fig1}
\end{figure}

Figure~\ref{fig1} shows how the expected value of $\AUC$ varies with
$n$. The dashed lines on either side of each curve are 95\% pointwise
confidence bands for the $\AUC$ estimate, and quantify the uncertainty
of the simulation study. The key feature is that, as predicted by
Theorem \ref{th3}, $\AUC$ changes very little with $n$, even when $n$ is small.
As expected, and as predicted by Theorems~\ref{th2}--\ref{th4}, $\AUC$ increases with
increasing $c_0$.

We also explored cases where the nonzero $\mu_j$'s took random values,
in particular where they were drawn randomly and uniformly from the
interval $[0,c_0 (20/n)^{1/2}]$. This case is more challenging, since
the genuine signals are now strictly smaller than in the previous
situation. Therefore, it comes as no surprise to learn that the $\AUC$
levels for each $c_0$ are reduced. However, the overall pattern of
stability with respect to $n$ is still evident, with very slightly more
variation than in the case of fixed $\mu_j$'s.

\subsection{\texorpdfstring{Influence of dependence on misranking performance}{In{f}luence of dependence on misranking performance}}\label{sec5.2}
Next, we discuss the effects of dependent noise $Z\ij$ in the model
at (\ref{eq3.0}). We take the noise to be a moving average of order 1,
that is, $Z_{ij} = \rho Z_{i,j-1} + (1-\rho^2)^{1/2} \ep_{i,j}$,
where the $\ep_{i,j}$'s are independent and normal N$(0,1)$. Thus,
$Z_{ij}$ and $Z_{ik}$ are correlated for all pairs $(j,k)$, with the
coefficient of correlation decaying exponentially fast in $|j-k|$. The
value of $c_0$ is fixed at 1.2, and nonzero $\mu_j$'s are chosen
uniformly in $[0,c_0 (20/n)^{1/2}]$. The values of $n$, $p$ and the
number of true signals are as in Section~\ref{sec5.1}.

Table~\ref{tab1} gives values of Monte Carlo approximations to the
means and standard deviations of $\AUC$ scores when the features for
which $\mu_j$ is nonzero are grouped together among the lowest values
of $j$, as in the discussion following Theorem \ref{th4}. The main observation
is that while mean $\AUC$ remains stable across the table, the
variability of $\AUC$ is much greater when strong dependence exists. For
instance, if $n=200$ then when $\rho=0.99$ the standard deviation of
$\AUC$ scores is ten times larger than when $\rho=0$. The presence of
dependence makes the problem significantly more difficult; it
effectively inserts an element of randomness into the process of
correctly ranking important features.

Table~\ref{tab2} shows results in the same setting, except that the
indices of features where $\mu_j$ is nonzero are distributed randomly
between 1 and $p$. There is again a high degree of stability, but
variability is comparatively less than that in Table~\ref{tab1},
consistent with the discussion in Section~\ref{sec3.3}. In particular, by
randomly distributing the indices of the nonzero $\mu_j$'s we
effectively reduce dependence among the features that are important,
and so, reflecting the results in Table~\ref{tab1}, the problem becomes
less statistically challenging.

\begin{table}
\caption{Mean and standard deviation of $\AUC$ scores for
simulation with correlated noise and grouped effects}\label{tab1}
\begin{tabular*}{\textwidth}{@{\extracolsep{\fill}}lllllllll@{}}
\hline
& \multicolumn{4}{l}{AUC means} & \multicolumn{4}{l@{}}{AUC
std dev.} \\[-6pt]
& \multicolumn{4}{l}{\hrulefill} & \multicolumn{4}{l@{}}{\hrulefill} \\
$\rho$ & $n=20$ & $n=50$ & $n=100$ & $n=200$ & $n=20$ & $n=50$ &
$n=100$ & $n=200$\\
\hline
$-0.99$& 0.725&0.698&0.706&0.709& 0.136&0.078&0.039&0.018 \\
$-0.75$& 0.704&0.717&0.715&0.712 & 0.064&0.022&0.012&0.006 \\
$-0.50$& 0.650&0.705&0.718&0.711 & 0.070&0.024&0.011&0.006 \\
$-0.25$& 0.702&0.725&0.706&0.716 & 0.062&0.025&0.012&0.007 \\
\phantom{$-$}0.00& 0.699&0.707&0.711&0.714 & 0.070&0.027&0.012&0.006 \\
\phantom{$-$}0.25& 0.687&0.714&0.707&0.712 & 0.080&0.032&0.015&0.007 \\
\phantom{$-$}0.50& 0.666&0.682&0.713&0.713 & 0.094&0.037&0.017&0.009 \\
\phantom{$-$}0.75& 0.715&0.710&0.718&0.719 & 0.101&0.047&0.025&0.013 \\
\phantom{$-$}0.99& 0.662&0.725&0.708&0.704 & 0.259&0.157&0.107&0.065 \\
\hline
\end{tabular*}
\end{table}

\begin{table}[b]
\caption{Mean and standard deviation of $\AUC$ scores for
simulation with correlated noise and randomised effects}\label{tab2}
\begin{tabular*}{\textwidth}{@{\extracolsep{\fill}}lllllllll@{}}
\hline
& \multicolumn{4}{l}{AUC means} & \multicolumn{4}{l@{}}{AUC
std dev.} \\[-6pt]
& \multicolumn{4}{l}{\hrulefill} & \multicolumn{4}{l@{}}{\hrulefill} \\
$\rho$ & $n=20$ & $n=50$ & $n=100$ & $n=200$ & $n=20$ & $n=50$ &
$n=100$ & $n=200$\\
\hline
$-0.99$& 0.712&0.712&0.713&0.715 & 0.139&0.058&0.032&0.016 \\
$-0.75$& 0.702&0.710&0.713&0.714 & 0.078&0.030&0.014&0.008 \\
$-0.50$& 0.705&0.710&0.711&0.714 & 0.076&0.032&0.015&0.008 \\
$-0.25$& 0.706&0.705&0.713&0.714 & 0.078&0.032&0.015&0.007 \\
\phantom{$-$}0.00& 0.696&0.709&0.714&0.712 & 0.077&0.030&0.015&0.007 \\
\phantom{$-$}0.25& 0.698&0.712&0.710&0.714 & 0.078&0.034&0.016&0.008 \\
\phantom{$-$}0.50& 0.703&0.713&0.714&0.714 & 0.081&0.030&0.016&0.007 \\
\phantom{$-$}0.75& 0.698&0.712&0.711&0.714 & 0.088&0.034&0.016&0.008 \\
\phantom{$-$}0.99& 0.741&0.721&0.712&0.710 & 0.195&0.091&0.053&0.024 \\
\hline
\end{tabular*}
\end{table}

\begin{figure}

\includegraphics{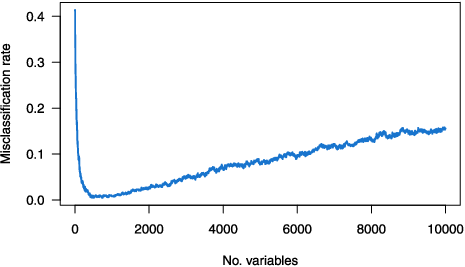}

\caption{Ideal prediction success for the example of
Section \protect\ref{sec5.3}.}\label{fig7.1}
\end{figure}

\subsection{Prediction in a large simulated problem}\label{sec5.3}
Here we present the analysis of a single simulated dataset,
demonstrating how our approach performs when the centroid classifier is
used. We take $p=10\mbox{, }000$ and $n=100$ (50 for each class). Ten percent
of the variables $X\ij$ (in the model at (\ref{eq3.0})) include a
nonzero signal. These $\mu_j$'s are drawn from the uniform distribution
on $[0,0.35]$, and the noise variables $Z\ij$ are independent $N(0,1)$.
This is a particularly difficult problem, since the signals are very
weak compared to the noise, and sample size is quite small.

Figure~\ref{fig7.1} shows the prediction performance of the centroid
classifier on a test set of 1000 replicates in the case of ``ideal
feature selection,'' where the 1000 features with nonzero signals are
selected first, in decreasing order of signal strength, followed by the
9000 features where the signal is not present. In particular, the
order is not chosen empirically. The minimum of the graph occurs at 449
features (out of a maximum of 1000), and corresponds to a
misclassification rate of only 0.5\%. The decrease in predictive
performance caused by less useful, or redundant, features is apparent
from the figure; the weaker genuine features actually hurt prediction
performance because they contain more noise than signal. Also of note
is the fact that a large number of features is needed to obtain good
prediction. For example, if attention is confined only to the strongest
50 features then the misclassification rate increases to 16\%.

\begin{figure}[b]

\includegraphics{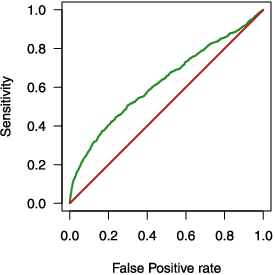}

\caption{ROC plot for feature ranking in the example of
Section \protect\ref{sec5.3}.}\label{fig7.2}
\end{figure}

For the same dataset, we undertook feature ranking based on the values
of $\hell_j$, defined at (\ref{eq2.3}). The result was the ROC chart in
Figure~\ref{fig7.2}. There the value of $k$ in the ranking $\hj
_1,\ldots
,\hj_k,\hj_{k+1},\ldots,\hj_p$, defined in the sentence below (\ref
{eq2.3}), is represented as $k/p$ on the horizontal axis, and the
vertical axis depicts the value of $\hell_{\hj_k}/\hell_{\hj_1}$, a
ratio of two negative numbers. The area under the empirical curve,
that is, $\AUC$, equals 0.626, meaning that a fraction $1-0.626\approx37\%$
of the paired scores correspond to a misranking. The ROC curve is
indexed by model size, with bottom left denoting an empty model and the
top right a full model. For a given model size, we can read off the
chart the corresponding sensitivity, or proportion of true important
features included, and the false positive rate, or proportion of
redundant features in the model. Ideally a model should have high
sensitivity and low false positive rate, and the chart indicates the
tradeoff between the two for various model sizes.

\begin{figure}

\includegraphics{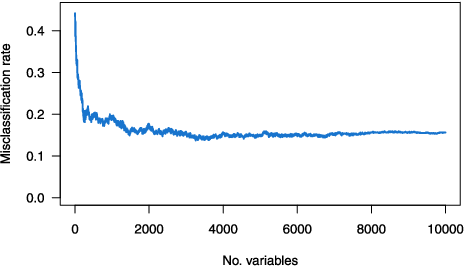}

\caption{Performance using feature ranking and centroid classifier
for the example of Section~\protect\ref{sec5.3}.}\label{fig7.3}
\end{figure}

Figure~\ref{fig7.3} shows how prediction accuracy varies with model
size. Performance is now clearly a long way from that represented in
Figure~\ref{fig7.1}, where the 1000 features with nonzero signals were
listed first in decreasing order of strength. The minimum misclassified
rate is now 13.7\%, and requires the use of 3258 of the 10,000
features. As discussed in the previous paragraph, every model size
corresponds to a position on the ROC plot in Figure~\ref{fig7.2}, in
this case $(0.31,0.51)$. Hence the optimal model found here contains 51\%
of the genuine features, and 31\% of the redundant ones.

We next explore stage (3) of the four-stage algorithm suggested in
Section~\ref{sec2.1}, addressing in turn each of the approaches (3a)--(3c)
discussed in Section~\ref{sec2.3}. We implement the threshold method, (3a), by
comparison with a randomised model where the observed classes $I_i$ are
scrambled and likelihoods recalculated, and employ the approach
suggested in the second paragraph of Section~\ref{sec4}; see (\ref{eq4.3}). In
particular, the threshold is chosen by computing the $100\a$th
percentile of scores for the scrambled data. Doing this for $\a=0.2$
corresponds to seeking a false positive rate of 0.2. In the numerical
example that we are considering here, this recovers a model with 2159
predictors and produces a test set misclassification rate of 16.8\%. A
model this size corresponds to the point (0.196, 0.40) on the ROC
chart. Notice that we have effectively targeted the false positive rate
of $\a=0.2$ via this approach.

\begin{figure}

\includegraphics{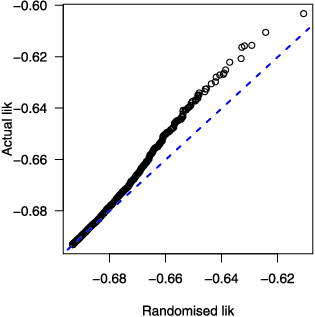}

\caption{Comparison of actual and randomised log-likelihoods in
the example of Section \protect\ref{sec5.3}.}\label{fig7.4}
\end{figure}

To provide an example of the change-point method, (3b), suggested in
Section~\ref{sec2.3} for choosing model size, we consider the ratio of the
sorted likelihoods from the original and scrambled rankings. These are
plotted in Figure~\ref{fig7.4}, along with a $45^\circ$ line. Starting
with the weakest features, we expect the ratio to remain near 1 until a
sizeable number of features that genuinely contain a positive signal
cause the ratio to shrink. For this purpose we can use the simple
change-point statistic for detecting a change in the mean (see
Chapter~2 of Cs\"org\H o and Horv\'ath \cite{CsoHor97}),
\[
T(t) = n^{-1/2} \bigl\{S(nt) - t S(n)\bigr\} ,
\]
where $S(k)$ equals the cumulative sum of the first $k$ ratios, and
$t\in(0,1)$ denotes the examined proportion of the dataset. This leads
to a model with 1760 features and a misclassification rate of 16.2\%,
comparable to that when using method (3a). This model size corresponds
to the point (0.16, 0.36) on the ROC plot.

Finally, in reference to the classifier-based approach (3c) suggested
in Section~\ref{sec2.3}, we note that the apparent error rate can be driven
quickly to zero without the actual error rate being reduced as much as
it is if we employ methods (3a) or (3b). For example, when using (3c)
in conjunction with the centroid classifier the ``best'' model, with
apparent error rate equal to zero, occurs when just 39 features are
selected; but the misclassification rate on the test set is 32.5\%,
almost twice that obtained for either of methods (3a) and~(3b).

\subsection{Results for microarray data examples}\label{sec5.4}
A challenge when using our methodology to analyse previously considered
real datasets is that the latter were possibly considered because they
illustrate cases where only a very small number of features determine
the class label. In particular, contrary to the concerns raised by
Goldstein~\cite{Gol09}, the number of influential components is quite
small. To simplify matters, we demonstrate here that likelihood based
ranking is a powerful tool for improving a wide variety of classifiers.
We make use of three well-known sets of microarray data. These relate,
respectively, to leukemia (Golub \etalc\cite{Goletal99}), colon cancer
(Alon \etalc \cite{Aloetal99}) and prostate cancer (Singh
\etalc~\cite{Singhetal02}) and have 7129, 2000 and 6033 components,
respectively. Dettling \cite{Det04} and Donoho and Jin~\cite{DonJin08}
discuss the performance of a variety of classifiers on these datasets,
using a two-thirds/one-third split of the data into training and test
samples. Their results are reported in Table~\ref{tabreal01}. Readers
may refer to the above papers for details on specific methods.

\begin{table}
\tablewidth=250pt
\caption{Percentage misclassification rate of different methods
on microarray datasets}\label{tabreal01}
\begin{tabular*}{250pt}{@{\extracolsep{\fill}}llll@{}}
\hline
Method & Leukemia & Colon & Prostate \\
\hline
Bagboo &4.08&16.10&\phantom{0}7.53\\
Boost &5.67&19.14&\phantom{0}8.71\\
RanFor &1.92&14.86&\phantom{0}9.00\\
SVM &1.83&15.05&\phantom{0}7.88\\
DLDA &2.92&12.86&14.18\\
KNN &3.83&16.38&10.59\\
PAM &3.55&13.53&\phantom{0}8.87\\
HCT &2.86&13.77&\phantom{0}9.47\\
\hline
\end{tabular*}
\end{table}

\begin{table}[b]
\caption{Performance of best methods on reduced datasets, using
likelihood based ranking} \label{tabrealb}
\begin{tabular*}{\textwidth}{@{\extracolsep{\fill}}lllllll@{}}
\hline
 & \multicolumn{2}{l}{Leukemia} & \multicolumn{2}{l}{Colon} &
\multicolumn{2}{l}{Prostate}\\[-6pt]
 & \multicolumn{2}{l}{\hrulefill} & \multicolumn{2}{l}{\hrulefill} &
\multicolumn{2}{l@{}}{\hrulefill}\\
Prop.& SVM & RanFor& RDA & RanFor& PAM & RanFor \\
\hline
0.0025&34.72&4.58&14.29&16.58&8.82&\phantom{0}8.00\\
0.005 &34.61&3.78&13.00&15.71&8.37&\phantom{0}7.67\\
0.01 &\phantom{0}6.86 &3.36&13.13&15.77&7.84&\phantom{0}7.73\\
0.025 &\phantom{0}1.89 &2.97&13.10&15.74&7.16&\phantom{0}8.20\\
0.05 &\phantom{0}1.67 &2.58&13.26&15.16&7.29&\phantom{0}8.55\\
0.10 &\phantom{0}1.58 &2.50&13.19&14.94&7.02&\phantom{0}8.90\\
0.15 &\phantom{0}1.67 &2.22&13.16&14.94&7.10&\phantom{0}9.14\\
0.25 &\phantom{0}1.64 &2.36&12.90&15.03&7.06&\phantom{0}9.49\\
0.375 &\phantom{0}1.61 &2.11&12.87&15.52&6.82&\phantom{0}9.67\\
0.50 &\phantom{0}1.58 &2.22&13.00&15.58&6.82&\phantom{0}9.90\\
0.75 &\phantom{0}1.53 &2.36&12.97&16.13&7.02&\phantom{0}9.80\\
1.00 &\phantom{0}1.61 &2.31&13.00&16.32&7.04&10.24\\
\hline
\end{tabular*}
\end{table}

To test the effectiveness of likelihood-based ranking, we chose the
best classification method and the random forest classifier (a
consistent performer) for each of the datasets. An extra step was added
to each cross-validation fold; the two-thirds training data was used to
rank features based on the likelihood score, and then only a proportion
of the top-ranked features were used to estimate the final model. The
results are presented in Table~\ref{tabrealb}. The last row of the
table shows results for the full dataset; they should in theory match
those in Table~\ref{tabreal01}, with differences attributable to tuning
approaches. We could not reproduce the accuracy reported for DLDA on
the colon dataset, and so used the next best method (PAM).

In each case, accuracy can be improved by reducing the model size. For
the best classifiers on each dataset, this effect was small but
noticeable; for the leukemia data, dimension was reduced by 25\% and
error by 5\%; for the colon dataset, dimension was reduced by 62.5\%
and error by 1\%; and for the prostate dataset, dimension was reduced
by 62.5\% and error by 3\%. For the random forest models, the results
were even more pronounced, with marked improvement in prediction and
significant dimension reduction. For the prostate dataset, the error
was reduced by 25\%, using just 0.005 of the available features in each
fold. This suggests that the likelihood based ranking method can
effectively control the sparsity of a model and potentially improve
model performance.

While firm conclusions are difficult here, we argue that this analysis
presents evidence for a large number of relatively weak effects
contributing to a model. Indeed, in all but one case we would prefer a
model size larger than the dozens, or fewer, used in many conventional
approaches to feature selection. Furthermore, our feature ranking
appears to be a useful means of determining the effective model size.

\subsection{Results for analysis of SNP data}\label{sec5.5}
We applied our methodology to SNP data collected to study multiple
sclerosis (ANZgene \cite{AusANZ09}; Wade \cite{Wad09}).
The original dataset consisted of 5031 subjects, which were collected
by two organisations: the Australian and New Zealand multiple sclerosis
genetic consortium (ANZgene) and the Wellcome Trust Multiple Sclerosis
Genetic Consortium 2. For data permission reasons, we restrict our
analysis to 3606 collected by ANZgene.
This part consists of 1618 case subjects (positive response) and 1988
control subjects (zero response).

For each of the subjects, we have 300,900 corresponding SNP
measurements, collected from different locations on the genome of each
patient. Each of these variables takes the values 0, 1 or~2, which we
treat as numeric. There were a small number of missing values, which we
addressed by setting them equal to the median value for that SNP.

\begin{figure}

\includegraphics{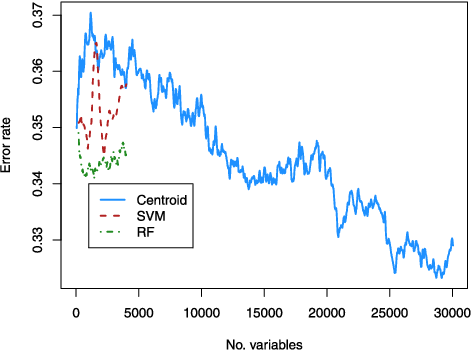}

\caption{Comparison of error rates on the multiple sclerosis SNP
data.}\label{SNP}
\end{figure}

The dataset was randomly divided into two-thirds training data and
one-third test data. Likelihood scores were calculated, and SNPs
ranked, based on their improvement over AIC. We then built three sets
of models:
\begin{itemize}
\item The centroid classifier using the top ranked $1,2,\ldots, 30\mbox{,}099$
variables. Note that this means that the final model was fitted using
10\% of the total variables;
\item A tuned SVM classifier using the top ranked $20, 40, 60, \ldots,
4000$ variables; and
\item A random forest classifier using the top ranked $20, 40, 60,
\ldots
, 4000$ variables.
\end{itemize}
The SVM and random forest models were not fitted beyond 4000 variables
because these were infeasible on desktop computers. For each model, the
classification error rate was measured on the test dataset. These
results are presented in Figure~\ref{SNP}.

We make two comments regarding these results. First, the trend in the
centroid classifier gives strong evidence of a large number of weak
effects; the performance continues to improve as more variables are
added, suggesting that there is useful information in the full top $10\%$
of the ranked variables. Second, the SVM and random forest approaches
performed better using a smaller number of SNPs, but worse when
compared to the large centroid models, again supporting the idea of
many weak effects. We would expect that adaptations of approaches like
SVM to large numbers of variables could also offer good performance.


\section*{Acknowledgements}
The authors would like to thank Melanie Bahlo, Lisa Melton, Terry
Speed, and Jim Stankovich for their help and generosity in sharing the
multiple sclerosis SNP data. We would also like to thank Li Liu for
help in preprocessing the SNP data.

%



\printhistory

\end{document}